\documentclass{article}
\usepackage[utf8]{inputenc}
\usepackage{amsmath,amssymb}
\usepackage{amsthm}
\usepackage{tikz}
\usetikzlibrary{arrows.meta}
\usepackage{parskip}
\usepackage{hyperref}
\usepackage[margin=1in]{geometry}
\usepackage{authblk} 
\usepackage{mathtools} 
\mathtoolsset{showonlyrefs,showmanualtags}
\usepackage{comment}

\usepackage{dsfont}
\usepackage{bbm}
\usepackage{mathrsfs}
\usepackage{graphicx}
\usepackage{multicol}
\usepackage{tikz}\usetikzlibrary{positioning,chains,fit,shapes,calc}

\mathtoolsset{showonlyrefs,showmanualtags}

\makeatletter
\def\thm@space@setup{%
  \thm@preskip=\parskip \thm@postskip=0pt
}
\makeatother

\hypersetup{
    colorlinks=true,
    linkcolor=blue,
    citecolor=darkgray
}

\usepackage{mathrsfs}
\usepackage{tikz}
\newtheorem{thm}{Theorem}[section]

\newtheorem{lemma}{Lemma}[section]

\newtheorem{prop}{Proposition}[section]
\newtheorem{remark}{Remark}
\theoremstyle{definition}

\usepackage[colorinlistoftodos]{todonotes}

\newcommand{\cH}{\mathcal{H}}

\newcommand{\EE}{\mathbb{E}}

\newcommand{\be}{\begin{equation}}
\newcommand{\ee}{\end{equation}}

\title{The length of the longest path in random acyclic orientations of a complete multipartite graph}

\author{Aiden Corcoran, Erik Lundberg, and Elisabeth Rodriguez}

\date{}

\begin{document}

\maketitle

\begin{abstract} 
We study the distribution of the length of the longest path in a uniformly random acyclic orientation of a complete multipartite graph. We derive an explicit analytic expression for the probability generating function of the longest path length, and, for an arbitrary fixed number of parts of equal size, we apply analytic combinatorics in several variables to obtain asymptotics (as the part size tends toward infinity) for its mean longest path length.
\end{abstract}



\section{Introduction}\label{sec:intro}

How long is the longest directed path in a uniformly random acyclic orientation of a fixed graph? This question, posed by Cameron \cite{CameronSlides}, asks for a probabilistic refinement of the enumeration of acyclic orientations. Khera and Lundberg \cite{KheraLundberg2025} addressed the case of a complete bipartite graph: they obtained a generating function for the distribution and, for two parts of equal size, proved a Gaussian limit law with explicit asymptotics for the mean and variance. They proposed extending the analysis to complete multipartite graphs \cite[Section~5]{KheraLundberg2025}. The present paper takes up that extension, obtaining a generating function providing the distribution for arbitrary part sizes and providing asymptotics for the mean when the parts have equal size.

An \emph{acyclic orientation} assigns a direction to each edge without creating a directed cycle. We measure the length of a directed path by its number of edges. Let $K_{n_1,\ldots,n_t}$ denote the complete multipartite graph with specified, labelled parts of sizes $n_1,\ldots,n_t$. Write $G_{n_1,\ldots,n_t}(\ell)$ for the number of its acyclic orientations whose longest directed path has length $\ell$, and put
\begin{equation}\label{eq:definition-F}
F(\boldsymbol z,u)
=\sum_{n_1,\ldots,n_t\geq0}\sum_{\ell\geq0}
G_{n_1,\ldots,n_t}(\ell)u^\ell
\prod_{j=1}^t\frac{z_j^{n_j}}{n_j!}.
\end{equation}
The empty graph contributes $1$; an edgeless nonempty graph has longest path length zero. Thus, as a generating function $F$ is exponential in the part-size variables and ordinary in the path-length variable.
The probability generating function (PGF) for a uniformly
sampled orientation is obtained by normalizing the
coefficient polynomial \cite[Chapter~III]{Flajolet}:
\begin{equation}\label{eq:PGF}
p_{n_1,\ldots,n_t}(u)
=\frac{[z_1^{n_1}\cdots z_t^{n_t}]F(\boldsymbol z,u)}
{[z_1^{n_1}\cdots z_t^{n_t}]F(\boldsymbol z,1)}.
\end{equation}

The evaluation $F(\boldsymbol z,1)$ gives the generating function for the total number of acyclic orientations (setting the marker variable $u=1$ unmarks the longest path length).  Stanley's identity expresses the number of acyclic orientations of a graph as an evaluation of its chromatic polynomial \cite{StanleyAO}. For complete bipartite graphs, Cameron and collaborators \cite{CameronPreprint} identified this number with a poly-Bernoulli number of negative index. These numbers were introduced by Kaneko \cite{Kaneko1997}; Brewbaker \cite{BrewbakerPaper} showed that they also enumerate lonesum binary matrices. Further combinatorial interpretations and bijections are developed by B\'enyi and Hajnal \cite{Benyi1,benyi2}. Building on earlier work of de Andrade, Lundberg, and
Nagle \cite{AndradeLundbergNagle} on the excedance set
statistic, Khera, Lundberg, and Melczer \cite{KLM}
obtained asymptotic enumeration results for lonesum
matrices, providing part of the analytic background
for the subsequent longest-path analysis in
\cite{KheraLundberg2025}.

More recently, Knuth established a bijection with so-called
\emph{parades} \cite{KnuthPreprint}, which will appear in the proof of Theorem \ref{thm:multipartiteIntro}.
Imagine $n_1$ girls and $n_2$ boys lining up for a parade,
with the members of each type labelled in increasing
order of age. The rule is that whenever two girls
march consecutively, the younger must come first,
and likewise for two consecutive boys. There is no
restriction on consecutive marchers of different types.
Thus a parade consists of alternating nonempty groups
of girls and boys, with the labels increasing within
each group.

To associate a graph with a parade, regard the marchers
as the vertices of $K_{n_1,n_2}$ and direct each edge from
the earlier marcher to the later one. The resulting
orientation is acyclic, since every directed edge
moves forward in the parade. Conversely, an acyclic
orientation recovers the parade by successively removing
its full set of sources and listing each such set in
increasing label order. This gives a bijection between
parades and acyclic orientations of $K_{n_1,n_2}$.
In \cite[Problem~8]{KnuthPreprint}, Knuth generalizes
this to parades with $t$ types of marchers, such as
girls, boys, and clowns when $t=3$, extending the
correspondence to complete $t$-partite graphs
$K_{n_1,\ldots,n_t}$.

Our first result provides an explicit analytic expression for the generating function $F(\boldsymbol z,u)$.

\begin{thm}\label{thm:multipartiteIntro}
For each integer $t\geq2$, the generating function \eqref{eq:definition-F} satisfies
\begin{equation}\label{eq:F-intro}
F(\boldsymbol z,u)
=\frac{1+(1-u)\displaystyle\sum_{j=1}^t
\frac{e^{z_j}-1}{1+u(e^{z_j}-1)}}
{1-u\displaystyle\sum_{j=1}^t
\frac{e^{z_j}-1}{1+u(e^{z_j}-1)}}.
\end{equation}
\end{thm}

For $t=2$, this reduces to the bipartite generating function of \cite{KheraLundberg2025},
\[
F(x,y,u)=\frac{e^{x+y}-(u-1)^2(e^x-1)(e^y-1)}
{1-u^2(e^x-1)(e^y-1)}.
\]
At $u=1$, it recovers Knuth's generating function
\begin{equation}\label{eq:Knuth-intro}
F(\boldsymbol z,1)=\frac1{\sum_{j=1}^t e^{-z_j} - (t-1)}
\end{equation}
for the number of acyclic orientations on a complete $t$-partite graph. 

Our second result gives asymptotics for the mean when the parts have equal size.

\begin{thm}\label{thm:mean}
Fix an integer $t\geq2$, and let $X_n$ be the longest directed path length in a uniformly random acyclic orientation of $K_{n,\ldots,n}$ with $t$ parts. Then
\begin{equation}\label{eq:mean}
\EE X_n = \frac{n}{\log(t/(t-1))} + O(1) ,\qquad n\longrightarrow\infty.
\end{equation}
\end{thm}

For $t=2$, the leading term is $n/\log2$, as in \cite{KheraLundberg2025}.
For large $t$, this asymptotic places the mean close to the maximal path length $tn-1$, since we have
\begin{equation}
\label{eq:logtasymp}
\frac{1}{\log(t/(t-1))}=t-\frac12+O(t^{-1}),
\qquad t\to\infty.
\end{equation}
More precisely, using \eqref{eq:logtasymp} and \eqref{eq:mean} we obtain
$$
\lim_{t\to\infty}\lim_{n\to\infty}
\frac{\EE X_n}{tn-1}=1.
$$

Let $B_{n_1,\ldots,n_t}=\sum_\ell G_{n_1,\ldots,n_t}(\ell)$ denote the total number of orientations. The proof of Theorem \ref{thm:mean} also yields the following result, giving asymptotics for the generalized parades in Knuth's question.

\begin{thm}\label{thm:KnuthAsymp}
Fix an integer $t\geq2$ and set $\rho=\log(t/(t-1))$. Then
\begin{equation}\label{eq:count}
B_{n,\ldots,n}
= (n!)^t\,
\frac{\dfrac{t}{t-1}(2\pi n)^{(1-t)/2}\rho^{-tn-1}}
{\sqrt{t(1-\rho)^{t-1}}}
\left(1+O(n^{-1})\right),
\qquad n\longrightarrow\infty.
\end{equation}
\end{thm}

\begin{remark}
This result was also obtained in the recent preprint \cite{Sun2026}.  We announced the result for the tripartite case at the 57th Southeastern International Conference on Combinatorics, Graph Theory, and Computing, in Spring 2026 \cite{CLR2026}.
\end{remark}

Section~\ref{sec:enumeration} proves Theorem~\ref{thm:multipartiteIntro} by decomposing orientations into successive source sets. Section~\ref{sec:asymptotics} proves Theorems~\ref{thm:mean} and~\ref{thm:KnuthAsymp} using smooth-point asymptotics from analytic combinatorics in several variables \cite{Melczer,PWM}.

\section{Acyclic orientations on complete multipartite graphs}\label{sec:enumeration}


The goal of this section is to prove Theorem \ref{thm:multipartiteIntro}, extending the bipartite generating function recalled in the introduction to the multipartite setting.
We follow the recursive approach that was used for the bipartite case in \cite[Remark~3]{KheraLundberg2025}, namely, pivoting on sets of sources. 
We begin with three parts to illustrate the derivation
of the system of equations. We then derive and solve
the system in the general case to obtain
Theorem~\ref{thm:multipartiteExplicit}.

\subsection{The tripartite case}\label{sec:tripartite}

The proof of Theorem~\ref{thm:multipartiteExplicit} uses a recursive decomposition related to Knuth's bijection between acyclic orientations and parades. A \emph{parade} is an ordering of labelled marchers in which labels increase within each maximal block of a single type. Under Knuth's bijection, these blocks correspond to the successive full sets of sources; see \cite[Section~7 and Problem~8]{KnuthPreprint} for the proof.

\begin{figure}[htbp]
\centering
\begin{tikzpicture}[
vertex/.style={circle,draw=black,fill=white,minimum size=6.2mm,inner sep=0pt,font=\scriptsize},
source/.style={draw=blue!70!black,fill=blue!15,line width=0.9pt},
edge/.style={-{Stealth[length=1.7mm,width=1.2mm]},draw=black!65,line width=0.45pt},
]
\begin{scope}[xshift=0cm]
\node[font=\small] at (1.175,1.95) {Initial orientation};
\node[vertex] (g1) at (0,1.15) {$g_1$};
\node[vertex,source] (g2) at (0,0) {$g_2$};
\node[vertex] (b1) at (2.35,1.15) {$b_1$};
\node[vertex] (b2) at (2.35,0) {$b_2$};
\node[vertex] (c1) at (0.58,-1.15) {$c_1$};
\node[vertex] (c2) at (1.77,-1.15) {$c_2$};
\draw[edge] (g1) -- (c1);
\draw[edge] (g1) -- (c2);
\draw[edge] (g2) -- (c1);
\draw[edge] (g2) -- (c2);
\draw[edge] (g2) -- (b1);
\draw[edge] (g2) -- (b2);
\draw[edge] (b1) -- (c1);
\draw[edge] (b1) -- (c2);
\draw[edge] (b1) -- (g1);
\draw[edge] (b2) -- (c1);
\draw[edge] (b2) -- (c2);
\draw[edge] (b2) -- (g1);
\node[font=\small,text=blue!70!black] at (1.175,-1.87) {$S_1=\{g_2\}$};
\end{scope}

\begin{scope}[xshift=3.85cm]
\node[font=\small] at (1.175,1.95) {After removing $S_1$};
\node[vertex] (g1) at (0,1.15) {$g_1$};
\node[vertex,source] (b1) at (2.35,1.15) {$b_1$};
\node[vertex,source] (b2) at (2.35,0) {$b_2$};
\node[vertex] (c1) at (0.58,-1.15) {$c_1$};
\node[vertex] (c2) at (1.77,-1.15) {$c_2$};
\draw[edge] (g1) -- (c1);
\draw[edge] (g1) -- (c2);
\draw[edge] (b1) -- (c1);
\draw[edge] (b1) -- (c2);
\draw[edge] (b1) -- (g1);
\draw[edge] (b2) -- (c1);
\draw[edge] (b2) -- (c2);
\draw[edge] (b2) -- (g1);
\node[font=\small,text=blue!70!black] at (1.175,-1.87) {$S_2=\{b_1,b_2\}$};
\end{scope}

\begin{scope}[xshift=7.7cm]
\node[font=\small] at (1.175,1.95) {After removing $S_2$};
\node[vertex,source] (g1) at (0,1.15) {$g_1$};
\node[vertex] (c1) at (0.58,-1.15) {$c_1$};
\node[vertex] (c2) at (1.77,-1.15) {$c_2$};
\draw[edge] (g1) -- (c1);
\draw[edge] (g1) -- (c2);
\node[font=\small,text=blue!70!black] at (1.175,-1.87) {$S_3=\{g_1\}$};
\end{scope}

\begin{scope}[xshift=11.55cm]
\node[font=\small] at (1.175,1.95) {After removing $S_3$};
\node[vertex,source] (c1) at (0.58,-1.15) {$c_1$};
\node[vertex,source] (c2) at (1.77,-1.15) {$c_2$};
\node[font=\small,text=blue!70!black] at (1.175,-1.87) {$S_4=\{c_1,c_2\}$};
\end{scope}

\end{tikzpicture}
\caption{Successive source removal in an acyclic orientation of $K_{2,2,2}$. The sources at each stage are highlighted. Removing $S_1,S_2,S_3,S_4$ in order exhausts the graph and gives the parade $g_2\mid b_1b_2\mid g_1\mid c_1c_2$, where bars separate maximal blocks of a single type. Its longest directed path has length $3$, for example $g_2\to b_1\to g_1\to c_1$.}
\label{fig:tripartite-source-removal}
\end{figure}
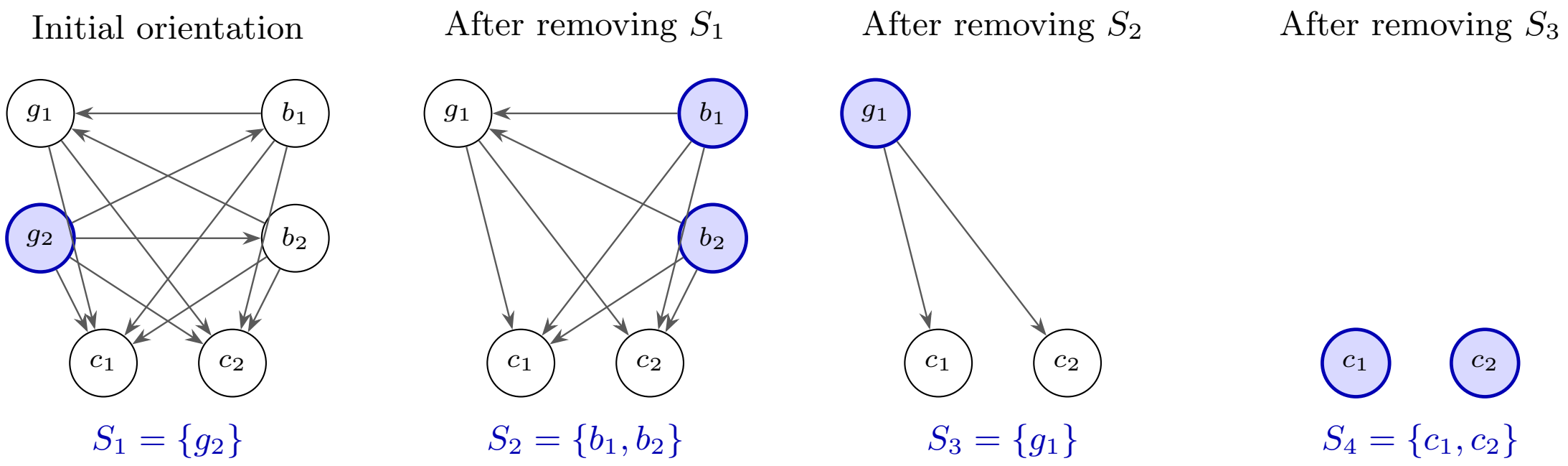

Let $F(x,y,z,u)$ denote the specialization of \eqref{eq:definition-F} to three parts. Following Knuth's narrative, call the vertices in the three parts girls, boys, and clowns, respectively. Let $G(x,y,z,u)$, $B(x,y,z,u)$, and $C(x,y,z,u)$ count nonempty acyclic orientations having a girl source, a boy source, and a clown source, in that order. These cases are disjoint, since no two vertices in different parts can both be sources. Thus, including the empty orientation,
\[
F(x,y,z,u)=1+G(x,y,z,u)+B(x,y,z,u)+C(x,y,z,u).
\]
The bijection constructs the parade by repeatedly removing each new set of sources, as illustrated in Figure~\ref{fig:tripartite-source-removal}. The figure shows an acyclic orientation of $K_{2,2,2}$ with parts $\{g_1,g_2\}$, $\{b_1,b_2\}$, and $\{c_1,c_2\}$. Its only source is $g_2$. After removing $g_2$, both $b_1$ and $b_2$ become sources; removing them makes $g_1$ a source. The final removal consists of the two isolated vertices $c_1,c_2$. Recording these successive sets in order, with increasing labels within each set, gives the parade $g_2\mid b_1b_2\mid g_1\mid c_1c_2$. The bars mark the four blocks, so the corresponding longest path has three edges.

Consider an orientation counted by $G(x,y,z,u)$ and decompose it according to the first step of the source removal process described above, namely, as its initial set of sources followed by an orientation on the graph that remains after source removal.  This leads to the following functional equation
\begin{equation}\label{eq:G-tripartite}
G(x,y,z,u)=(e^x-1)\Bigl(1+u\bigl(B(x,y,z,u)+C(x,y,z,u)\bigr)\Bigr).
\end{equation}
Indeed, the initial block is a nonempty labelled set of girls, giving the factor $e^x-1$. Removing it leaves either the empty orientation, contributing $1$, or an orientation whose sources are boys or clowns, counted by $B(x,y,z,u)+C(x,y,z,u)$. In the latter case, restoring the initial block increases the longest path length by one, giving the factor $u$.

Conversely, given an orientation whose sources are
boys or clowns, we may introduce a nonempty labelled
set of girls, with all their edges directed outward. The new girls are
exactly the sources of the resulting orientation.

By symmetry, interchanging the three types of marchers and their corresponding variables gives the system
\begin{equation}\label{eq:tripartite-system}
\left\{
\begin{aligned}
G(x,y,z,u)&=(e^x-1)\Bigl(1+u\bigl(B(x,y,z,u)+C(x,y,z,u)\bigr)\Bigr),\\
B(x,y,z,u)&=(e^y-1)\Bigl(1+u\bigl(G(x,y,z,u)+C(x,y,z,u)\bigr)\Bigr),\\
C(x,y,z,u)&=(e^z-1)\Bigl(1+u\bigl(G(x,y,z,u)+B(x,y,z,u)\bigr)\Bigr).
\end{aligned}
\right.
\end{equation}
We solve this system after generalizing it to the multipartite case below.

\subsection{The general multipartite case}

The above result has the following generalization to the multipartite setting.

\begin{lemma}\label{lem:source-system}
For $j=1,\ldots,t$, let $G_j(\boldsymbol z,u)$ be the generating function for nonempty acyclic orientations of complete $t$-partite graphs whose sources lie in part $j$. Then
\begin{equation}\label{eq:source-decomposition}
F(\boldsymbol z,u)=1+\sum_{j=1}^tG_j(\boldsymbol z,u),
\end{equation}
and the functions $G_j$ satisfy the linear system
\begin{equation}\label{eq:source-system}
G_j(\boldsymbol z,u)
=(e^{z_j}-1)\left(1+u\sum_{i\ne j}G_i(\boldsymbol z,u)\right),
\qquad j=1,\ldots,t.
\end{equation}
\end{lemma}

\begin{proof}
Every nonempty acyclic orientation has a source, and its sources all lie in one part (since vertices in different parts are adjacent). This proves \eqref{eq:source-decomposition}.

Fix an orientation counted by $G_j$, and let $S$ be its full set of sources. If removing $S$ leaves the empty graph, the orientation is supported on a nonempty set in part $j$ and has no edges. This gives the empty-remainder contribution $1$ inside the parentheses in \eqref{eq:source-system}.
Otherwise, removing $S$ leaves an acyclic orientation whose sources belong to a part $i\ne j$.

Moreover, removing $S$ decreases the longest path length by exactly one. Indeed, any longest path in the graph that remains (after removing $S$) starts at a source, which lies in part $i\ne j$, and can therefore be extended by a vertex of $S$. Conversely, any directed path in the original orientation uses at most one vertex of $S$, necessarily as its first vertex. In both cases, choosing the nonempty labelled set $S$ contributes the common factor $e^{z_j}-1$. The remainder contributes either $1$ or $u\sum_{i\ne j}G_i(\boldsymbol z,u)$, so the labelled product construction \cite{Flajolet} gives \eqref{eq:source-system}.
\end{proof}

Solving the system \eqref{eq:source-system} leads to Theorem \ref{thm:multipartiteIntro}, restated here for convenience.

\begin{thm}\label{thm:multipartiteExplicit}
The generating function $F(\boldsymbol z,u)$ has the explicit analytic expression
\begin{equation}\label{eq:F}
F(\boldsymbol z,u)
=
\frac{
1 + (1-u)\displaystyle\sum_{j=1}^t \frac{e^{z_j}-1}{1+u(e^{z_j}-1)}
}{
1 - u\displaystyle\sum_{j=1}^t \frac{e^{z_j}-1}{1+u(e^{z_j}-1)}
}.
\end{equation}
\end{thm}

\begin{proof}
The system is linear in $G_j$.
For ease of notation, we write
$a_j = e^{z_j} - 1, \quad G = \sum_{j=1}^t G_j$.
Adding $u a_j G_j$ to both sides of the $j$th equation, our system becomes
$$
G_j (1 + u a_j) = a_j (1 + u G), \quad j = 1,2,...,t.
$$
Dividing by $(1+u a_j)$ and summing over $j$ gives a single equation involving $G$:
$$
G = (1 + u G)\sum_{j=1}^t \frac{a_j}{1 + u a_j}.
$$
Solving for $G$ gives
$$
G =
\frac{\sum_{j=1}^t \frac{a_j}{1+u a_j}}
{1 - u \sum_{j=1}^t \frac{a_j}{1+u a_j}},
$$
and hence
$$
F(\boldsymbol z,u)
=
1 + G
=
\frac{1 + (1-u)\sum_{j=1}^t \frac{a_j}{1+u a_j}}
{1 - u \sum_{j=1}^t \frac{a_j}{1+u a_j}},
$$
as desired.
\end{proof}

\section{Asymptotic enumeration and mean longest path length}\label{sec:asymptotics}

\subsection{Preliminary results from ACSV}\label{sec:ACSV-prelim}

Let $F(\boldsymbol z)=G(\boldsymbol z)/H(\boldsymbol z)^p$, where $p\geq1$ is an integer,
$G,H$ are entire functions, $H(\boldsymbol0)\neq0$, and common factors
have been canceled. Write
\[
F(\boldsymbol z)=\sum_{\boldsymbol\alpha\in\mathbb N^d}f_{\boldsymbol\alpha} \boldsymbol z^{\boldsymbol\alpha},
\qquad \boldsymbol z=(z_1,\ldots,z_d),
\]
and let
\[
\mathcal V=\{\boldsymbol z\in\mathbb C^d:H(\boldsymbol z)=0\}
\]
denote its singular variety. Throughout, we consider points of
$\mathcal V$ with nonzero coordinates.

A point $\boldsymbol w\in\mathcal V$ is called \emph{minimal} if there is no
$\boldsymbol z\in\mathcal V$ such that $|z_j|<|w_j|$ for every $j$.
It is called \emph{strictly minimal} if
\[
\mathcal V\cap
\{\boldsymbol z\in\mathbb C^d:|z_j|\leq |w_j|\text{ for every }j\}
=\{\boldsymbol w\}.
\]
The point $\boldsymbol w$ is \emph{smooth} if $\nabla H(\boldsymbol w)\neq0$.

Fix $\boldsymbol r=(r_1,\ldots,r_d)\in\mathbb Q_{>0}^d$.
A smooth point $\boldsymbol w\in\mathcal V$ is called a \emph{critical point in the direction $\boldsymbol r$}
if
\be\label{eq:crit}
\bigl(w_1H_{z_1}(\boldsymbol w),\ldots,w_dH_{z_d}(\boldsymbol w)\bigr)
=\lambda \boldsymbol r
\ee
for some $\lambda\neq0$.
In all coefficient asymptotics below, $n\to\infty$ through
positive integers for which $n\boldsymbol r\in\mathbb N^d$.

Suppose $\boldsymbol w$ is a smooth critical point and
$H_{z_d}(\boldsymbol w)\neq0$. By the implicit function theorem,
$\mathcal V$ is locally given by
\[
z_d=g(\widehat{\boldsymbol z}),
\qquad
\widehat{\boldsymbol z}=(z_1,\ldots,z_{d-1}).
\]
Set $\widehat{\boldsymbol w}=(w_1,\ldots,w_{d-1})$ and
$\widehat{\boldsymbol r}=(r_1,\ldots,r_{d-1})$, and define
\[
\phi(\boldsymbol\theta)=
\log\frac{
g(w_1e^{i\theta_1},\ldots,w_{d-1}e^{i\theta_{d-1}})
}{g(\widehat{\boldsymbol w})}
+i\frac{\widehat{\boldsymbol r}\cdot\boldsymbol\theta}{r_d},
\]
using the logarithm with $\phi(\boldsymbol0)=0$.
The critical point equations \eqref{eq:crit} give $\nabla\phi(\boldsymbol0)=0$.
Let $\cH$ denote the Hessian matrix of $\phi$ at $\boldsymbol\theta=\boldsymbol0$. The critical point $\boldsymbol w$ is
\emph{nondegenerate} if $\det\mathcal H\neq0$.

We will need the following result referred to as smooth-point
asymptotics \cite[Theorem 5.3]{Melczer}.

\begin{thm}[Smooth asymptotics for poles]\label{thm:ACSV-poles}
With $F(\boldsymbol z)=G(\boldsymbol z)/H(\boldsymbol z)^p$ as above, where $p\geq1$ is an integer,
suppose $\boldsymbol w\in\mathcal V\cap\mathbb R_{>0}^d$ is a strictly
minimal smooth critical point in direction $\boldsymbol r$, with
\[
G(\boldsymbol w)\neq0,\qquad H_{z_d}(\boldsymbol w)\neq0,
\]
and suppose $\mathcal H$ is real and nonsingular.
Then
\[
f_{n\boldsymbol r}
=
\frac{G(\boldsymbol w)}
{(p-1)!\bigl(-w_dH_{z_d}(\boldsymbol w)\bigr)^p}
\frac{
r_d^{p-1}\boldsymbol w^{-n\boldsymbol r}n^{p-1}(2\pi n)^{(1-d)/2}
}{
\sqrt{\det(r_d\mathcal H)}
}
\left(1+O(n^{-1})\right).
\]
\end{thm}

Note: Minimality implies that $\Re\phi$ has a local minimum
at the origin. Hence, the real Hessian $\mathcal H$
is positive semidefinite, and the assumption that $\cH$ is nonsingular then makes it positive definite.

A power series $P(\boldsymbol z)=\sum_{\boldsymbol\alpha\in\mathbb N^d}p_{\boldsymbol\alpha} \boldsymbol z^{\boldsymbol\alpha}$ is called \textit{aperiodic} if every element of $\mathbb{Z}^d$ can be written as an integer linear combination of the exponent vectors  
$\{\boldsymbol\alpha:p_{\boldsymbol\alpha}\neq0\}$ that appear in $P$.

The following proposition (\cite[Prop. 5.5]{Melczer}, \cite[Prop. 3.17]{PemWilTwenty}) will be used to verify the strictly minimal condition.

\begin{prop}\label{prop:aperiodicity}
Suppose $H$ is analytic on $\mathbb C^d$ and $H=1-P$,
where $P$ is an aperiodic power series with nonnegative
coefficients and $P(\boldsymbol0)<1$. Then every minimal point of
$\mathcal V=\{\boldsymbol z:H(\boldsymbol z)=0\}$ with nonzero coordinates is
strictly minimal and has positive real coordinates.
\end{prop}

\begin{remark}\label{rmk:multiply}
Note that the definition of strict minimality depends only on the variety $\mathcal V$.  Hence, before applying the proposition, it is permissible to multiply $H$ by a nonvanishing entire function.
\end{remark}

\subsection{Proof of Theorems \ref{thm:mean} and \ref{thm:KnuthAsymp}}

Recall from \eqref{eq:PGF} that the PGF $p_n(u):=p_{n,\ldots,n}(u)$ for the longest path length (where we are specializing to the case of parts of equal size $n$) can be expressed as a ratio of coefficient extractions.
The mean is then obtained by
evaluating the derivative at $u=1$, i.e., we have \cite[Chapter~III]{Flajolet}
\begin{equation}\label{eq:ratio}
\EE X_n
=p_n'(1) = \frac{[\boldsymbol z^{n\boldsymbol1}]F_u(\boldsymbol z,1)}
{[\boldsymbol z^{n\boldsymbol1}]F(\boldsymbol z,1)},
\end{equation}
where $\boldsymbol z^{n\boldsymbol1}:=z_1^n z_2^n \cdots z_t^n$.

We derive our asymptotic results by applying Theorem \ref{thm:ACSV-poles} to obtain asymptotics for the numerator and denominator of \eqref{eq:ratio}.

First let us evaluate $F(\boldsymbol z,u)$ and its partial derivative
$F_u(\boldsymbol z,u)$ at $u=1$.
Differentiating \eqref{eq:F} with respect to $u$, we obtain
\[
F_u(\boldsymbol z,u)
=
\frac{
\left(\displaystyle\sum_{j=1}^t
\frac{e^{z_j}-1}{1+u(e^{z_j}-1)}\right)^2
-\displaystyle\sum_{j=1}^t
\frac{(e^{z_j}-1)^2}{(1+u(e^{z_j}-1))^2}
}{
\left(1-u\displaystyle\sum_{j=1}^t
\frac{e^{z_j}-1}{1+u(e^{z_j}-1)}\right)^2
}.
\]
Setting $u=1$ gives
\begin{equation}\label{eq:mean-generating-functions}
F(\boldsymbol z,1)=\frac1{H(\boldsymbol z)},
\qquad
F_u(\boldsymbol z,1)=\frac{A(\boldsymbol z)}{H(\boldsymbol z)^2},
\end{equation}
where
\begin{equation}\label{eq:H-A}
H(\boldsymbol z)=1-t+\sum_{j=1}^t e^{-z_j},
\qquad
A(\boldsymbol z)=\left(\sum_{j=1}^t(1-e^{-z_j})\right)^2
-\sum_{j=1}^t(1-e^{-z_j})^2.
\end{equation}

The equal part size condition corresponds to the diagonal direction $\boldsymbol r=\boldsymbol1$.
Hence, the critical point equations \eqref{eq:crit} become
\[
H(\boldsymbol w)=0,
\qquad
w_1e^{-w_1}=\cdots=w_te^{-w_t}.
\]
By symmetry, we seek a positive solution with
$w_1=\cdots=w_t=\rho$. Applying this and recalling \eqref{eq:H-A}, the equation $H(\boldsymbol w)=0$
then reduces to
\[
1-t+te^{-\rho}=0,
\]
giving
\[
\rho=\log\frac{t}{t-1},
\qquad
\boldsymbol w=(\rho,\ldots,\rho)
\]
as the candidate for the contributing singularity.
The next lemma verifies that this critical point
is smooth and strictly minimal.

\begin{lemma}\label{lem:critical}
The point $\boldsymbol w$ is a strictly minimal smooth
critical point in the diagonal direction for both functions
in \eqref{eq:mean-generating-functions}.
\end{lemma}

\begin{proof}
Following the observation in Remark \ref{rmk:multiply}, before we apply Proposition~\ref{prop:aperiodicity}, we multiply $H$ by an entire nonvanishing function so that the result has the form $1-P$, where $P$ has nonnegative
coefficients:
\[
e^{z_1+\cdots+z_t}H(\boldsymbol z)
=\sum_{j=1}^t e^{\sum_{i\ne j}z_i}
-(t-1)e^{z_1+\cdots+z_t}
=1-P(\boldsymbol z).
\]
Then $P(\boldsymbol0)=0$, and the coefficient of
$\boldsymbol z^{\boldsymbol\alpha}$ in $P$ for
$\boldsymbol\alpha\ne\boldsymbol0$ is
$(k-1)/\boldsymbol\alpha!$, where $k$ is the number of
positive coordinates of $\boldsymbol\alpha$.
Thus $P$ has nonnegative coefficients, with support
consisting of all exponent vectors having at least two
positive coordinates. This support generates $\mathbb Z^t$, since
\[
\boldsymbol e_i
=(2\boldsymbol e_i+\boldsymbol e_j)
-(\boldsymbol e_i+\boldsymbol e_j)
\qquad(j\ne i).
\]
Hence $P$ is aperiodic.

Since $H(\boldsymbol w)=0$, we have $P(\boldsymbol w)=1$.
Nonnegativity of the coefficients of $P$ gives
\[
|P(\boldsymbol z)|
\le P(|z_1|,\ldots,|z_t|)
<P(\boldsymbol w)=1, \quad |z_j|<w_j, \quad j=1,2,...,t.
\]
This proves $\boldsymbol w$ is minimal.
Proposition~\ref{prop:aperiodicity} therefore gives strict
minimality.

Finally, $H_{z_j}(\boldsymbol w)=-(t-1)/t\ne0$ for every $j$,
so $\boldsymbol w$ is smooth and critical in the diagonal
direction. Both numerators in \eqref{eq:mean-generating-functions} are nonzero there, since
$A(\boldsymbol w)=(t-1)/t$.
\end{proof}

\begin{lemma}\label{lem:hessian}
The phase Hessian $\cH$ at $\boldsymbol w$ in the diagonal
direction is
\begin{equation}\label{eq:mean-hessian}
\mathcal H=(1-\rho)(I_{t-1}+J_{t-1}),
\end{equation}
where
$I_{t-1}$ is the identity matrix
and $J_{t-1}$ is the all-ones matrix.
In particular, its determinant is given by
\begin{equation}
\det\mathcal H=t(1-\rho)^{t-1}.
\end{equation}
\end{lemma}

\begin{proof}
Locally write $H=0$ as $z_t=g(\widehat{\boldsymbol z})$.
Implicit differentiation of
\[
e^{-g(\widehat{\boldsymbol z})}
+\sum_{j=1}^{t-1}e^{-z_j}=t-1
\]
at $\widehat{\boldsymbol w}=(\rho,\ldots,\rho)$ gives
\[
g(\widehat{\boldsymbol w})=\rho,\qquad
\frac{\partial g}{\partial z_j}(\widehat{\boldsymbol w})=-1,
\qquad
\frac{\partial^2 g}{\partial z_j\partial z_k}
(\widehat{\boldsymbol w})=1+\delta_{jk},
\]
where $\delta_{jk}$ is the Kronecker delta.
Setting $h=\log g$, we obtain
\[
\frac{\partial h}{\partial z_j}(\widehat{\boldsymbol w})
=-\frac1\rho,
\qquad
\frac{\partial^2 h}{\partial z_j\partial z_k}
(\widehat{\boldsymbol w})
=\frac{1+\delta_{jk}}{\rho}-\frac1{\rho^2}.
\]
In the diagonal direction, the phase is
\[
\phi(\boldsymbol\theta)
=h(\rho e^{i\theta_1},\ldots,\rho e^{i\theta_{t-1}})
-\log\rho+i\sum_{j=1}^{t-1}\theta_j.
\]
Differentiating twice with respect to the angular
variables and evaluating at $\boldsymbol\theta=\boldsymbol0$
therefore gives
\[
\mathcal H_{jk}
=\frac{\partial^2\phi}{\partial\theta_j\partial\theta_k}
(\boldsymbol0)
=-\rho^2
\frac{\partial^2 h}{\partial z_j\partial z_k}
(\widehat{\boldsymbol w})
-\delta_{jk}\rho
\frac{\partial h}{\partial z_j}(\widehat{\boldsymbol w})
=(1-\rho)(1+\delta_{jk}).
\]
Thus, $\mathcal H=(1-\rho)(I_{t-1}+J_{t-1})$.
Applying a formula \cite[Lemma~1.1]{DingZhou} for the determinant of a rank-one update of an invertible matrix gives
\[
\det(I_{t-1}+J_{t-1})=1+(t-1)=t.
\]
Hence, $\det\cH=t(1-\rho)^{t-1}$.
\end{proof}

\begin{proof}[Proof of Theorem~\ref{thm:KnuthAsymp}]
Apply Theorem~\ref{thm:ACSV-poles} to
$F(\boldsymbol z,1)=1/H(\boldsymbol z)$ with
$p=1$, $d=t$, and $\boldsymbol r=\boldsymbol1$.
Lemma~\ref{lem:critical} verifies
its hypotheses. Since $H_{z_t}(\boldsymbol w)=-(t-1)/t$, and Lemma~\ref{lem:hessian} gives $\det\cH=t(1-\rho)^{t-1}$,
we obtain
\begin{equation}\label{eq:denominator-asymptotic}
[\boldsymbol z^{n\boldsymbol1}]F(\boldsymbol z,1)
=
\frac{\dfrac{t}{t-1}\rho^{-tn-1}(2\pi n)^{(1-t)/2}}
{\sqrt{t(1-\rho)^{t-1}}}
\left(1+O(n^{-1})\right).
\end{equation}
Multiplying by $(n!)^t$ proves \eqref{eq:count}.
\end{proof}

\begin{proof}[Proof of Theorem~\ref{thm:mean}]
Recalling \eqref{eq:ratio} we have
\[
\EE X_n
=\frac{[\boldsymbol z^{n\boldsymbol1}]F_u(\boldsymbol z,1)}
{[\boldsymbol z^{n\boldsymbol1}]F(\boldsymbol z,1)}.
\]
Applying Theorem~\ref{thm:ACSV-poles} to
$F_u(\boldsymbol z,1)=A(\boldsymbol z)/H(\boldsymbol z)^2$
with $p=2$ and $\boldsymbol r=\boldsymbol1$, we obtain
\[
[\boldsymbol z^{n\boldsymbol1}]F_u(\boldsymbol z,1)
=
\frac{nA(\boldsymbol w)}
{\rho^2 H_{z_t}(\boldsymbol w)^2}
\frac{\rho^{-tn}(2\pi n)^{(1-t)/2}}
{\sqrt{t(1-\rho)^{t-1}}}
\left(1+O(n^{-1})\right).
\]
Dividing by \eqref{eq:denominator-asymptotic}, the
exponential and Gaussian factors cancel, giving
\[
\EE X_n
=-\frac{nA(\boldsymbol w)}
{\rho H_{z_t}(\boldsymbol w)}+O(1).
\]
At $\boldsymbol w$, each $1-e^{-w_j}=1/t$, so
\[
A(\boldsymbol w)=1-\frac1t
=\frac{t-1}{t}=-H_{z_t}(\boldsymbol w).
\]
Consequently $\EE X_n=n/\rho+O(1)$, proving
\eqref{eq:mean}.
\end{proof}

\section{Concluding remarks}

\begin{figure}[htbp]
\centering
\begin{tabular}{@{}ccc@{}}
\includegraphics[width=.18\textwidth]{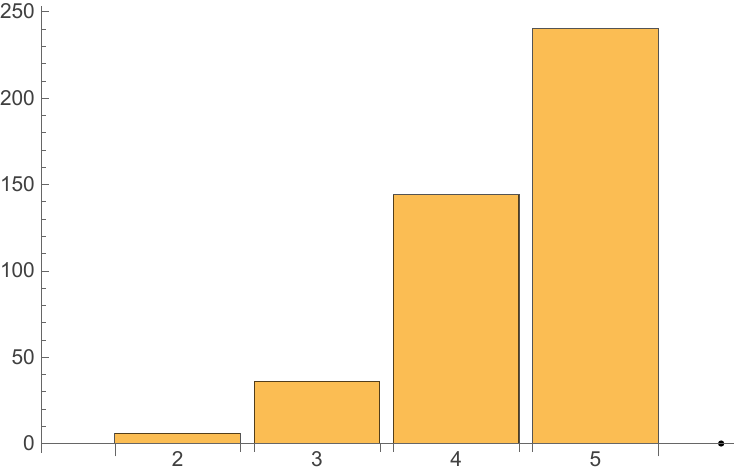}
&
\includegraphics[width=.2\textwidth]{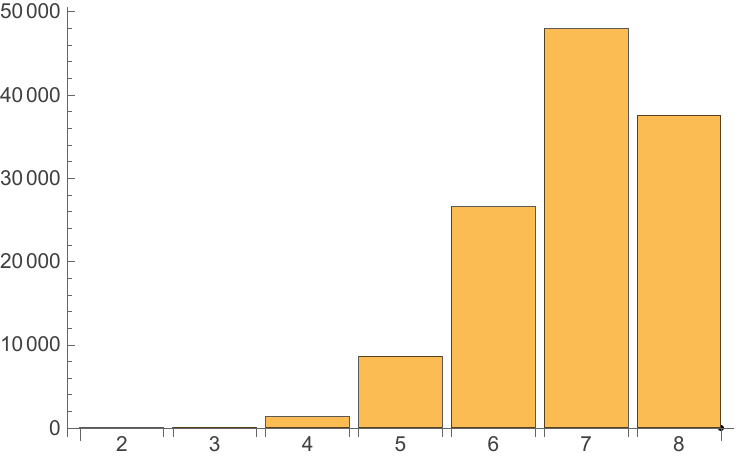}
&
\includegraphics[width=.2\textwidth]{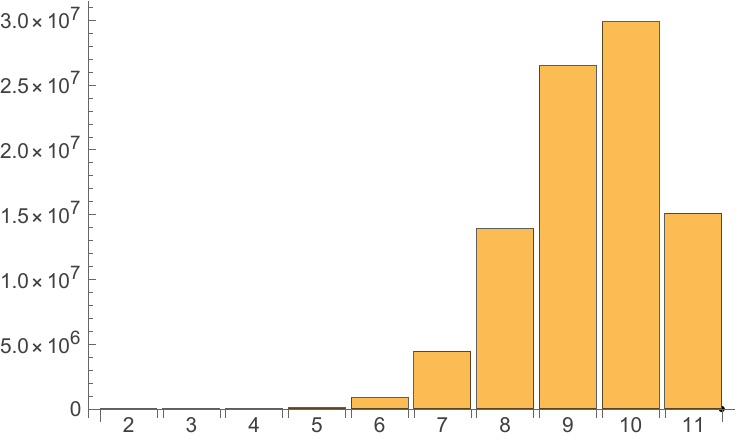}
\\
$n=2$ & $n=3$ & $n=4$
\\[1ex]
\includegraphics[width=.27\textwidth]{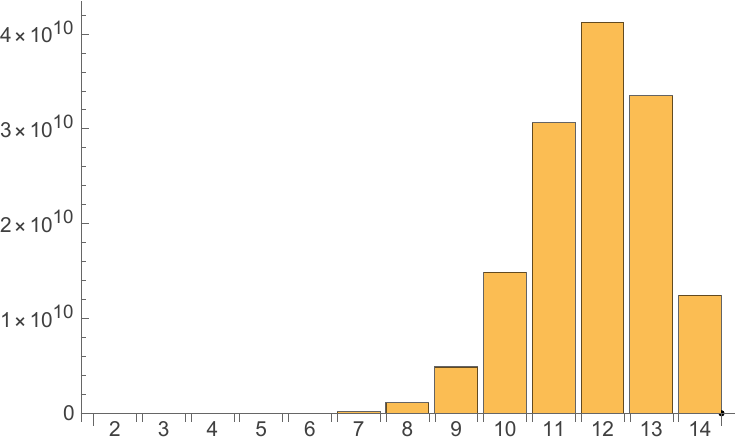}
&
\includegraphics[width=.27\textwidth]{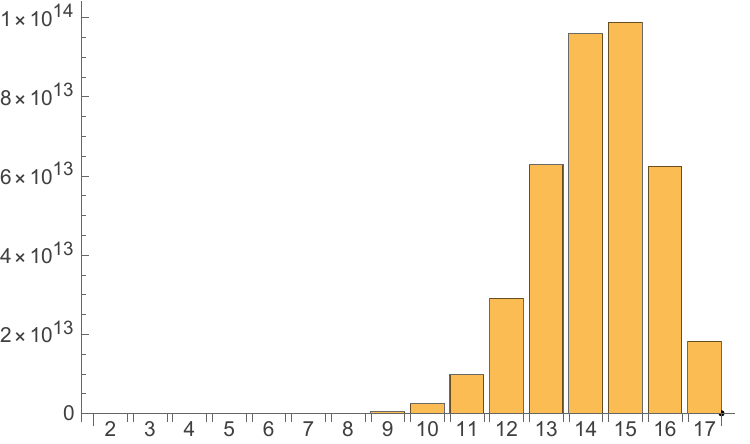}
&
\includegraphics[width=.283\textwidth]{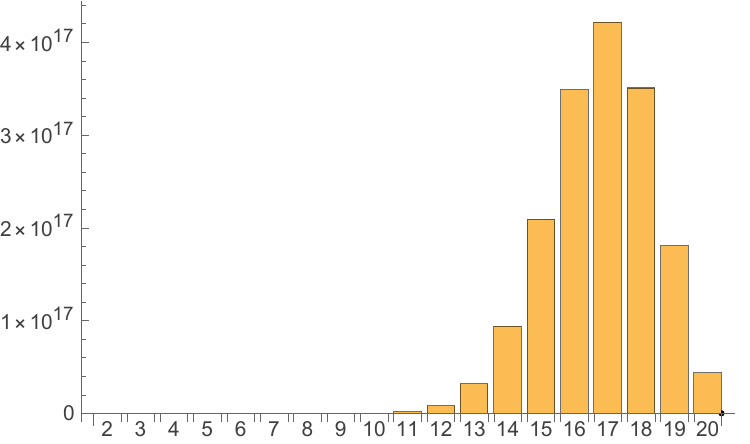}
\\
$n=5$ & $n=6$ & $n=7$
\\[1ex]
\includegraphics[width=.3\textwidth]{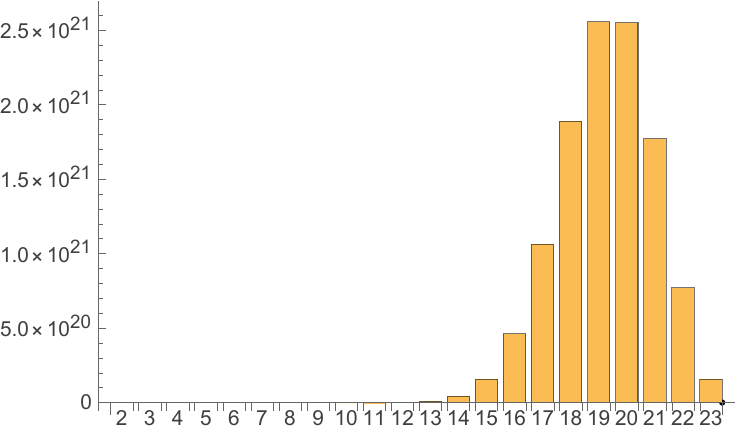}
&
\includegraphics[width=.31\textwidth]{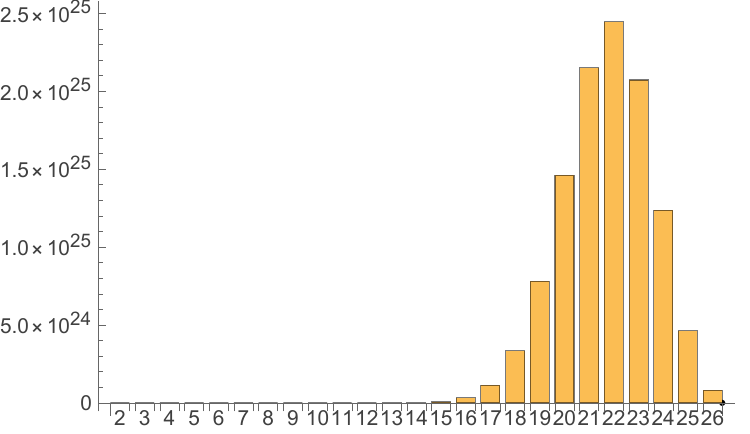}
&
\includegraphics[width=.31\textwidth]{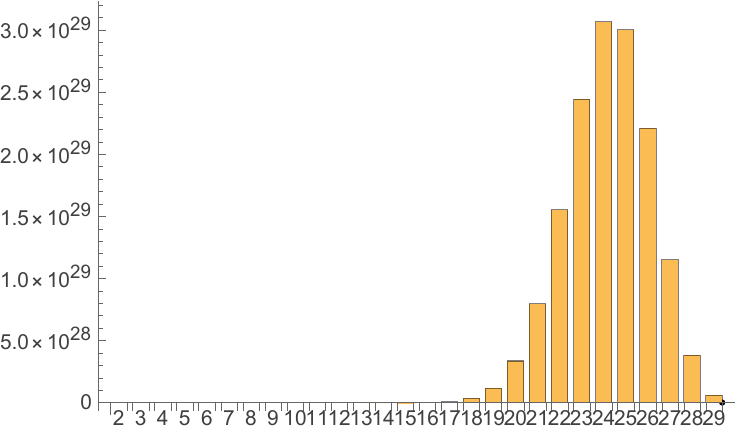}
\\
$n=8$ & $n=9$ & $n=10$
\end{tabular}
\caption{Distributions of the longest path length in
uniformly sampled acyclic orientations of $K_{n,n,n}$
for $n=2,\ldots,10$, computed by exact
coefficient extraction from \eqref{eq:F}. These distributions
provide numerical evidence for a Gaussian limit after
centering and scaling.}
\label{fig:gaussian-evidence}
\end{figure}

For complete bipartite graphs with equal part sizes tending
to infinity, the longest path length is asymptotically
Gaussian \cite{KheraLundberg2025}. It is natural to ask
whether this extends to complete $t$-partite graphs for
every fixed $t\ge3$. Writing $X_{t,n}$ for the longest
path length in a uniformly sampled acyclic orientation
with $t$ parts of size $n$, we conjecture
\[
\frac{X_{t,n}-\EE X_{t,n}}
{\sqrt{\operatorname{Var}(X_{t,n})}}
\xrightarrow{d}\mathcal N(0,1)
\qquad(n\to\infty).
\]
The generating function \eqref{eq:F} provides a starting
point for this question. Beyond the asymptotic mean determined in Theorem \ref{thm:mean}, one would need to determine the asymptotic variance
and control the coefficient asymptotics uniformly as $u$ varies near $u=1$.
Figure~\ref{fig:gaussian-evidence} provides numerical
evidence for this conjecture in the tripartite case.
A further question is whether a Gaussian limit holds
when the part sizes grow in fixed, unequal proportions.

In the alternative regime,
where the number of parts $t$ increases with the part size $n$ fixed,
it seems more natural to expect Poisson rather than Gaussian
behavior after suitable recentering of the longest path length.
More generally, one may let both the part size and the
number of parts increase at related rates and look for a possible phase transition in the limiting longest path length distribution as the relative growth varies.

\noindent\textbf{Acknowledgments and AI tools disclosure.} 
All mathematical results and their proofs were conceived and written by the authors.
We used
OpenAI GPT-6.1 to proofread the paper and suggest improvements in the exposition. 
The second-named author acknowledges support from Simons Foundation grant 712397. 

\bibliographystyle{abbrv}
\bibliography{acyclic}

\bigskip
\noindent
Department of Mathematics and Statistics\\
Florida Atlantic University\\
Boca Raton, FL 33431, USA

\medskip
\noindent
Aiden Corcoran: \texttt{acorcoran2023@fau.edu}\\
Erik Lundberg: \texttt{elundber@fau.edu}\\
Elisabeth Rodriguez: \texttt{elisabethrod2024@fau.edu}

\end{document}